\documentclass [twoside,10pt]{article}
\usepackage{graphicx}
\usepackage{amssymb}
\usepackage{amsthm}
\newtheorem{thm}{Theorem}

\newtheorem{lem}{Lemma}
\newtheorem{pro}{Proposition}
\newtheorem{cla}{Claim}
\newcommand{\bd}{{\rm bd}}
\newcommand{\diam}{{\rm diam}}

\newcommand{\width}{{\rm width}}
\newcommand{\dist}{{\rm dist}}

\newcommand{\sh}{{\rm sinh \,}}
\newcommand{\ch}{{\rm cosh \,}}
\newcommand{\arsh}{{\rm arcsinh \,}}
\newcommand{\arch}{{\rm arccosh \,}}

\title{\bf Width of Convex Bodies in Hyperbolic Space}

\date{}

\begin{document}

\maketitle

\thispagestyle{empty}

\vskip-1.2cm

\noindent{\author{Marek Lassak}}

\pagestyle{myheadings} \markboth{\centerline {Marek Lassak}}{\centerline {Width of convex bodies in hyperbolic space}}

\baselineskip 12.5pt

\maketitle
\vskip 0.5cm

\noindent 
{\bf Abstract}.
For every hyperplane $H$ supporting a convex body $C$ in the hyperbolic space $\mathbb{H}^d$ we define the width of $C$ determined by $H$ as the distance between $H$ and a most distant ultraparallel hyperplane supporting $C$. 
We define bodies of constant width in $\mathbb{H}^d$ in the standard way as bodies whose all widths are equal.
We show that every body of constant width is strictly convex.
The minimum width of $C$ over all supporting $H$ is called the thickness $\Delta (C)$ of $C$.
A convex body $R \subset \mathbb{H}^d$ is said to be reduced if $\Delta (Z) < \Delta (R)$ for every convex body $Z$ properly contained in $R$.
We show that regular tetrahedra in $\mathbb{H}^3$ are not reduced.
Similarly as in the Euclidean and spherical spaces, we introduce complete bodies and bodies of constant diameter also in $\mathbb{H}^d$. 
We show that every body of constant width $\delta$ is a body of constant diameter $\delta$ and a complete body of diameter $\delta$.
Moreover, the two last conditions are equivalent.

\baselineskip 16.5pt



\medskip

\date{}

\maketitle

\section{Introduction} 

Here we present a sketch of the content of this paper. 

In Section 2 we introduce some notions and present four lemmas.

Section 3 is devoted to the notion of width of a convex body $C$ in the hyperbolic space $\mathbb{H}^d$.
For any hyperplane $H$ supporting $C$ we define ${\rm width}_H(C)$ of $C$ determined by $H$ as the distance between $H$ and any farthest ultraparallel hyperplane supporting $C$. 
Proposition 1 shows that ${\rm width}_H (C)$ equals to the maximum distance between $H$ and a point of $C$, and Proposition 2 claims that ${\rm width}_H (C)$ equals to the distance between $H$ and the nearest equidistant surface $E$ to $H$ such that $C$ is a subset of the equidistant strip bounded by $H$ and $E$. 

Section 4 is about the diameter and thickness.
Theorem 1 says that the maximum width of $C$ equals to the diameter of $C$.
The thickness $\Delta (C)$ of a convex body $C \in \mathbb{H}^d$ is defined as the infimum of ${\rm width}_H (C)$ over all hyperplanes $H$ supporting $C$.
In Theorem 2 we assume that $\width_H (C) = \Delta (C)$ for a hyperplane $H$ supporting $C$ and that there exists exactly one farthest point $j$ of $C$ from $H$.
The thesis is that the projection of $j$ onto $H$ belongs to $C \cap H$.
Finally, we find a formula for the width of the regular simplex determined by the hyperplane containing its facet.

Section 5 concerns bodies of constant width and reduced bodies.
We define the notion of a body of constant width as a body whose all widths are equal.
In Proposition 3 we show that every body of constant width is strictly convex.
By a reduced body we mean a convex body $R \subset \mathbb{H}^d$ such that $\Delta (Z) < \Delta (R)$ for every convex body $Z \subset R$ different from $R$. 
Proposition 4 says that every body of constant width is a reduced body.
In Remark 2 we explain why every regular tetrahedron in $\mathbb{H}^3$ is not reduced.

In Section 6 we consider complete bodies and bodies of constant diameter in $\mathbb{H}^d$.
Theorem 3 says that every body of constant width $\delta$ is a body of constant diameter $\delta$, and that a body is of constant diameter $\delta$ if and only if it is a complete body of diameter $\delta$. 
We conjecture that the three kinds of bodies in $\mathbb{H}^d$ coincide. 

In our considerations it is convenient to work with the hyperboloid model of $\mathbb{H}^d$, so about the model on the upper sheet $x_{d+1}= \sqrt{x_1^2 + \dots + x_d^2 +1}$ of the two-sheeted hyperboloid.
This model permits a reasonable comparison of the obtained results with the analogous ones in the Euclidean space $\mathbb{E}^d$ and the spherical space $\mathbb{S}^d$.
Our figures show the orthogonal look to the sheet from the above.

\section{A few notions and lemmas}

The geodesic between two different points $a, b \in \mathbb{H}^d$ is called the {\bf segment} $ab$.  
By a {\bf ball of radius} $\rho$ we mean the set of points of $\mathbb{H}^d$ which are at a distance at most $\rho$ from a fixed point called the {\bf center} of this ball.
If $d=2$, the ball is called a {\bf disk}, and its boundary a {\bf circle}.
We say that a set $C \subset \mathbb{H}^d$ is {\bf convex} if together with every two points $a,b$ it contains the whole segment $ab$.  
By a {\bf convex body} in $\mathbb{H}^d$ we mean a closed bounded convex set with non-empty interior. 
Of course, the intersection of any family of convex sets is also convex. 
Thus for every set $Q \subset \mathbb{H}^d$ there exists the unique smallest convex set containing $Q$. 
It is called {\bf the convex hull of} $Q$.

If a hyperplane $G$ of $\mathbb{H}^d$ has a common point with a convex body $C \subset \mathbb{H}^d$ and if its intersection with the interior of $C$ is empty, we say that $G$ {\bf supports} $C$.
If at every boundary point of $C$ exactly one hyperplane supports $C$, the body is said to be {\bf smooth}.

Let $A, B \subset \mathbb{H}^d$ be two non-empty sets.
The symbol ${\rm dist} (A, B)$ denotes the infimum of distances $|ab|$ over all $a \in A$ and $b \in B$.
We call it the distance between $A$ and $B$.
If $A$ is a one point set $\{c\}$, we simply write ${\rm dist} (c, B)$.

We omit an easy proof of the following lemma.

\begin{lem}
For the projection $h$ of a point $g$ onto a hyperplane $H$ we have $\dist(g, H) = |gh|$.
Moreover, $|gk| > |gh|$ for every $k \in H$ different from $h$.
\end{lem}

Recall that two hyperplanes $H$ and $J$ in $\mathbb{H}^d$ are said to be {\bf ultraparallel} if they do not intersect in a real or imaginary point. 
It is well known that they have exactly one common orthogonal straight line and that the distance of its intersections with $H$ and $J$ is ${\rm dist} (H, J)$. 
The convex hull of $H \cup J$ is called an {\bf ultraparallel strip}.

\begin{lem} 
Let $H$ and $J$ be ultraparallel hyperplanes and let $p$ be a point in this half-space bounded by $H$ which does not contain $J$. 
Then ${\rm dist} (p, H) < {\rm dist} (p, J)$.
\end{lem}

\begin{proof}  
Clearly, $|ph| = {\rm dist} (p, H)$ for the projection $h$ of $p$ onto $H$ and $|pj| = {\rm dist} (p, J)$ for the projection $j$ of $p$ onto $J$.
Denote by $g$ the point of intersection of $pj$ with $H$.
By Lemma 1 we have $|ph| \leq |pg|$.
Hence ${\rm dist} (p, H) = |ph| \leq |pg| < |pj| = {\rm dist} (p, J)$.
\end{proof}

\begin{lem} 
Assume that $H$ supports a convex body $C$ and let $j$ be a farthest point of $C$ from $H$. 
Denote by $h$ the projection of $j$ onto $H$.
Let $J$ be the hyperplane through $j$ which is orthogonal to $jh$.
Then $J$ supports $C$ at $j$.
Moreover, $J \cap C = \{j\}$.
\end{lem}

\begin{proof} 
Imagine that $J$ does not support $C$ (see Fig. 1).
Then a point $k \in J$ belongs to the interior of $C$.
By Lemma 1 for our $J$ in part of $H$ there, from $k \not = j$ we obtain $\dist(k,H,) > \dist(j,H)$.
Let $g$ be the projection of $k$ onto $H$.
There is a point $m \in \bd (C)$ such that $k \in gm$.
Clearly $m \not = k$ and thus $\dist(m,H) > \dist(k,H)$.
Hence $\dist(m,H) > \dist(j,H)$. 
A contradiction with the choice of $j$.
Consequently, $J$ supports $C$.
Moreover, $J \cap C$ cannot contain a point different from $j$, since such a point would be in a larger distance from $H$ than $j$.
\end{proof}

\vskip0.1cm
\begin{figure}[htbp]        
\hskip0.6cm
\includegraphics[width=5.58cm]{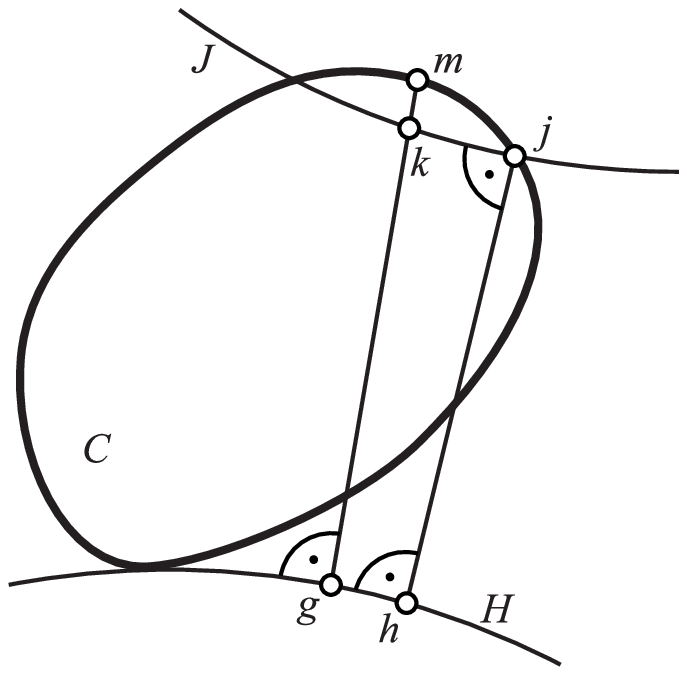} \hskip0.8cm 
\includegraphics[width=6.84 cm]{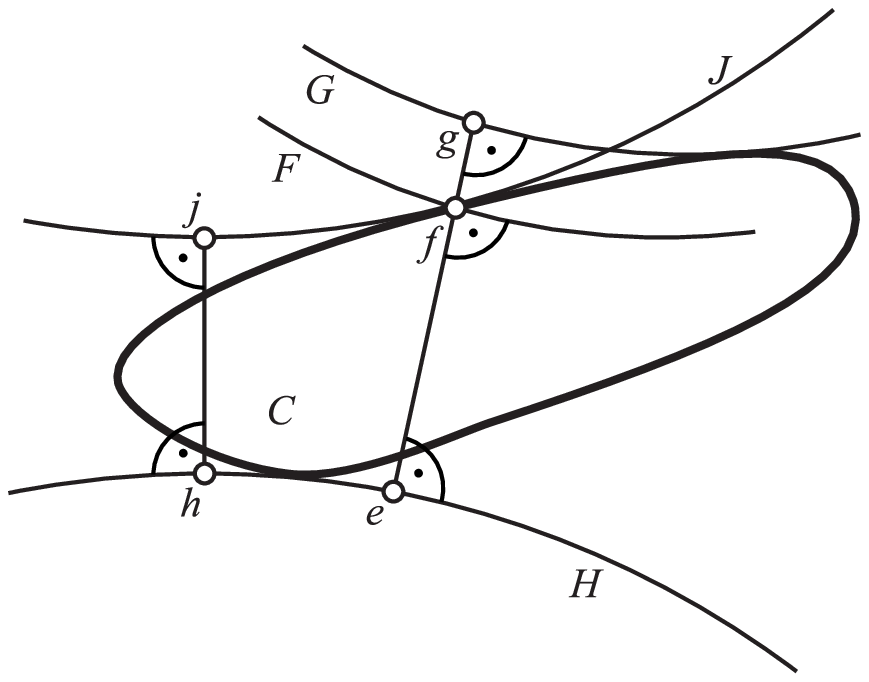} 

{\bf Fig. 1.} Illustration to the proof of Lemma 3 \hskip 0.15cm 
{\bf Fig. 2.} Illustration to the proof of Lemma 4   

\end{figure}

\begin{lem} 
Assume that a hyperplane $H$ supports a convex body $C$ and let $J$ be 
a supporting ultraparallel hyperplane of $C$ whose distance from $H$ is maximal. 
Then $J \cap C$ consists of exactly one point.
\end{lem}

\begin{proof}
There are $h \in H$ and $j \in J$ such that $hj$ is orthogonal to $H$ and $J$.
Clearly, $\dist (H, J) = |hj|$.

Imagine that $J \cap C$ has more than one point.
Then there is at least one point $f \in J \cap C$ different from $j$ (see Fig. 2).
Denote by $e$ the projection of $f$ onto $H$.
Clearly, $|ef| > |hj|$.

Let $F$ be the hyperplane through $f$ orthogonal to the segment $ef$.
By Lemma 1 we have $\dist (F, H) = |ef| > |hj| = \dist (J, H)$.

If $F$ supports $C$, then it is a farthest ultraparallel hyperplane to $H$ supporting $C$ than $J$.
This contradicts the assumption that $J$ is a farthest ultraparallel hyperplane to $H$ supporting $C$.

If $F$ does not supports $C$, then $F$ passes through the interior of $C$.
There exists a supporting hyperplane $G$ of $C$ orthogonal to the straight line containing $ef$ such that $F$ is a subset of the interior of the strip between $H$ and $G$.
By $g$ denote the point of intersection of $G$ with this line.

Of course, $g \not = f$.
Thus by  $f \in eg$ we obtain $|ef| < |eg|$.
So $\dist (G, H) = |ge|$ and $\dist (F, H) = |fe|$ imply $\dist (G, H) > \dist (F, H)$. 
Thus from the earlier established inequality $\dist (F, H) > \dist (J, H)$ we get $\dist (G, H) > \dist (J, H)$. 
This contradicts the description of $J$ as the most distant ultraparallel to $H$ hyperplane supporting $C$.
The obtained contradiction implies that $J \cap C$ is a one-point set.
\end{proof}

Assume that a hyperplane $H$ supports a convex body $C \in \mathbb{H}^d$ and $J$ is a farthest from $H$ ultraparallel hyperplane supporting $C$.
Then a most distant from $J$ ultraparallel hyperplane $K$ supporting $C$ does not have to be $H$ (differently than in $\mathbb{E}^d$).
A simple example showing this is the regular simplex and the hyperplane $H$ containing a facet of~it.

Let $p$ be a point out of a hyperplane $H$.
Recall that the {\bf equidistant surface to $H$ through $p$} is the set of all points in the distance $\dist(p,H)$ from $H$ which are in the half-space bounded by $H$ and containing $p$.

\section{Width of a convex body} 

Let $H$ be a hyperplane supporting a convex body $C \subset \mathbb{H}^d$. 
We define the {\bf width of $C$ determined by $H$} as the distance between $H$ and any farthest ultraparallel hyperplane supporting $C$ (by compactness arguments there exists at least one such a most distant one, sometimes there are a finite or even infinitely many of them). 
The symbol ${\rm width}_H (C)$ denotes this width of $C$ determined by $H$.

\begin{pro} 
Let $C \subset \mathbb{H}^d$ be a convex body and $H$ be any supporting hyperplane of $C$.
Then $\width_H (C)$ equals to the maximum distance between $H$ and a point of $C$. 
\end{pro}

\begin{proof}
Since $C$ is compact, there exists a point $j \in C$ in the maximum distance from $H$.
Denote by $h$ the projection of $j$ onto $H$.
Provide the hyperplane $J$ through $j$ which is orthogonal to $hj$.
Clearly, $H$ and $J$ are ultraparallel.
By Lemma 3 we conclude that $J$ supports $C$ at $j$ (see Fig.~3).

\vskip0.1cm
\begin{figure}[htbp]        
\hskip3.9cm
\includegraphics[width=6.48cm]{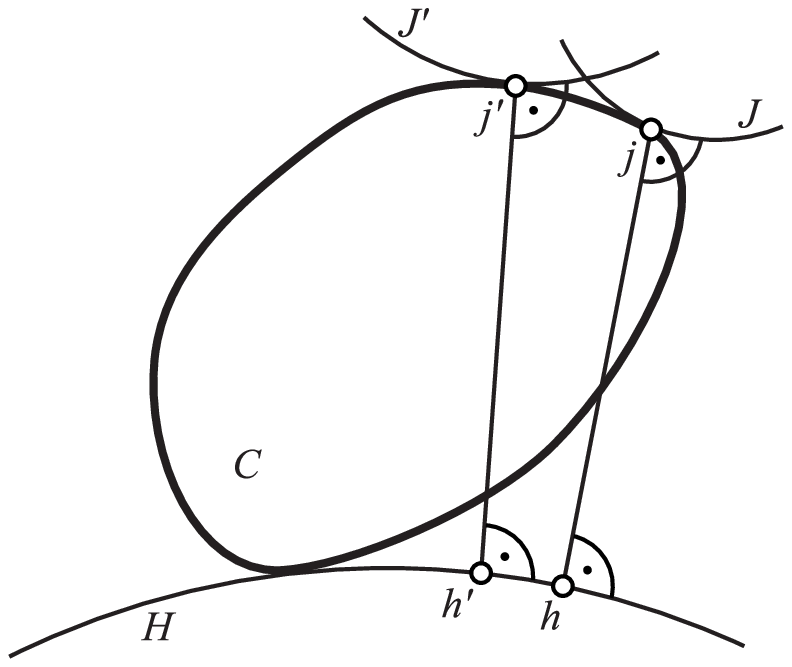} 

\vskip0.5cm
\hskip 3.1cm {\bf Fig. 3.} Illustration to the proof of Proposition 1 

\end{figure}

Let us show that $J$ is a farthest ultraparallel hyperplane to $H$ supporting $C$.
Imagine the opposite. 
Then there is a more distant from $H$ ultraparallel hyperplane $J'$ supporting $C$.
Denote by $j'$ its closest to $H$ point of support. 
Its projection onto $H$ is denoted by $h'$.
By $\dist(H, J') > \dist (H, J)$ we have $|h'j'| > |hj|$.
This contradicts the choice of $j$.

We have shown that $J$ is a farthest ultraparallel hyperplane to $H$ supporting $C$.
Thus by the definition of ${\rm width}_H (C)$ we have ${\rm width}_H (C) = \dist (H, J)$.
Since $h \in H$, $j \in J$ and $hj$ is orthogonal to both these ultraparallel hyperplanes, we have $\dist (H, J) = \dist (H, j)$.
Hence ${\rm width}_H (C) = \dist (H, j)$, which is our thesis.
\end{proof}

In connection with this proposition, pay attention that in all the geometries $\mathbb{E}^d$, $\mathbb{S}^d$ and $\mathbb{H}^d$, we may define the ${\rm width}_H (C)$ as the maximum distance between $H$ and a point of $C$.
Recall here that ${\rm width}_H (C)$ in $\mathbb{E}^d$ is the width of $C$ in the direction orthogonal to $H$, and ${\rm width}_H (C)$ for $\mathbb{S}^d$ is defined in \cite{[L15]}.

We omit an easy proof of the following proposition.

\begin{pro} 
Let $C \subset \mathbb{H}^d$ be a convex body and $H$ any supporting hyperplane of $C$.
Then ${\rm width}_H (C)$ equals to the distance between $H$ and the nearest equidistant surface $E$ to $H$ such that $C$ is a subset of the equidistant strip being the convex hull of $H \cup E$. 
\end{pro}

Each of the above propositions may be regarded as a different definition of ${\rm width}_H (C)$. 
 
Let us add that a few different notions of width of a convex body $C \subset \mathbb{H}^d$ are known.
For a comparison of them see the paper \cite{[Ho]} by Horv\'ath.

\section{Diameter and thickness of a convex body}

We put the material on the diameter and the thickness in one section since they are of analogous nature.
Just we see this from the below Theorem 1 characterizing the diameter of a convex body $C \subset \mathbb{H}^d$ and the below definition of the thickness of $C$: they are the maximum and the minimum values of $\width_H(C)$ over all hyperplanes $H$ supporting $C$. 

By the {\bf diameter} ${\rm diam}(A)$ of a set $A \subset \mathbb{H}^d$ we mean the supremum of the distances between pairs of points of $A$.
Clearly, if $A$ a convex body, by its compactness, ${\rm diam}(A)$ is realized for at least one pair of points of $A$.

\begin{thm} 
For every convex body $C \subset \mathbb{H}^d$ we have 
\vskip -0.4cm

$$\max \{ {\rm width}_H (C); H {\rm \ is \ a \ supporting \ hyperplane \ of} \ C \} = {\rm  diam}(C).$$
\end{thm}

\begin{proof}
Denote the left side of this equality by $\max \width (C)$.

By the compactness of $C$ there exists a supporting hyperplane $K$ of $C$ such that $\width_K(C) = \max \width (C)$.
By Proposition 1 
the value $\width_K(C)$ equals to the maximum distance of a point $c \in C$ from $K$.
For the projection $k$ of $c$ onto $K$ we have $|kc| = \width_K(C)$.
Since $K$ supports $C$, there exists a point $f \in K \cap C$.
Of course $|cf| \geq |ck|$.
Hence $\diam(C) \geq \max \width (C)$.

There are $h, j$ in $\bd(C)$ such that $|hj| = \diam(C)$.
Provide the hyperplane $H$ through $h$ orthogonal to $hj$.
Observe that $H$ supports $C$ (still if not, then there exists a point $y \in C$ which is in the strictly opposite side of $H$ than $j$, which means that $|yj| > |hj|$ in contradiction to $|hj| = \diam (C)$).
Since $H$ supports $C$ and the point $j$ is in $C$, by Proposition 1 we get $\width_H (C) \geq |hj|$.
Hence $\max \width (C) \geq |hj|$.
This and $|hj| = \diam (C)$ lead to $\diam(C) \leq \max \width (C)$.

From the above two paragraphs we obtain the thesis of our theorem.
\end{proof}

\begin{cla} 
Assume that ${\rm  diam} (C) = |ab|$ for some points $a, b$ of a convex body $C \subset \mathbb{H}^d$. 
Denote by $S$ the ultraparallel strip whose both bounding hyperplanes are orthogonal to the segment $ab$. 
We have $C \subset S$. 
\end{cla}

\begin{proof}
Denote by $A$ the hyperplane through $a$ bounding $S$ and by $B$ the hyperplane through $b$ bounding $S$.
Imagine that $C \not \subset S$.
Then a point $c \in C$ is out of $S$.
It is strictly separated from $S$ by $A$ or $B$.
For instance, let it be separated by $A$.
Denote by $a'$ the intersection of $bc$ with $A$.
Of course, $a' \in C$.
From the right triangle $baa'$ we conclude that $|a'b| > |ab|$.
This contradicts ${\rm  diam} (C) = |ab|$.
Consequently, $C \subset S$. 
\end{proof}

By {\bf the thickness $\Delta (C)$ of} a convex body $C \subset \mathbb{H}^d$ we mean the infimum of ${\rm width}_H (C)$ over all hyperplanes $H$ supporting $C$.
By compactness arguments, this infimum is realized, so $\Delta (C)$ is the minimum of the numbers ${\rm width}_H (C)$.

\vskip0.2cm
\noindent
{\bf Remark 1.}
In $\mathbb{E}^d$ the following statement holds true. 
{\it Let $C$ be a convex body.
Let $H$ be a supporting hyperplane of $C$ such that $\width_H(C) = \Delta (C)$.
Denote by $j$ any farthest point of $C$ from $H$. 
Then the projection of $j$ onto $H$ belongs to $H \cap C$.}
It follows from the property formulated by Eggleston at the bottom of page 77 of \cite{[Eg]}. 
An analogous statement holds true in $\mathbb{S}^d$ as a consequence of Claim 2 of \cite{[L15]} (see also Corollary 2.13 of \cite{[L22]}).
Submitting the first version of this paper, the author expected that such a statement holds also in $\mathbb{H}^d$.
But it is not true as it results from the following example shown by the referee.

\vskip0.2cm
\noindent
{\bf Example 1.} Let $H$ be a line and $E$ be an equidistant curve to $H$ in a distance $\rho$. 
Take a point $a \in H$ and the orthogonal line $\ell$ to $H$ through $a$.
Let $b$ and $c$ be points of $E$ symmetric with respect to $\ell$ such that the angle $\angle bac$ is right or obtuse. 
Define $C$ as the convex hull of the set consisting of the point $a$ and all the points of $E$ between $b$ and $c$.
Clearly, $\width_H (C) = \rho = \Delta(C)$.
The projections of all the points of $E$ between $bc$ besides its middle onto $H$ do not belong to $H \cap C$.

\vskip 0.2cm
Here is a weaker form of the (not true in $\mathbb{H}^d$) statement from Remark 1.

\begin{thm} 
Let $C \subset \mathbb{H}^d$ be a convex body and let $H$ be a supporting hyperplane of $C$ such that $\width_H (C) = \Delta (C)$.
Assume that there exists a unique most distant point $j \in C$ from $H$.
Then the projection of $j$ onto $H$ belongs to $H \cap C$.
\end{thm} 

\begin{proof}
Imagine the opposite to the thesis, this is that the projection $h$ of $j$ onto $H$ does not belong to $H \cap C$.

Since $H \cap C$ is convex and compact, there exists exactly one point $g \in H \cap C$ in the minimum distance from $h$. 
The above uniqueness of $g$ results from the fact that $g$ is the intersection of the $(d-1)$-dimensional ball $B \subset H$ centered at $h$ which touches $H \cap C$.
We do not loose the generality assuming that the vector $hg$ is directed into the positive orientation on the straight line containing the segment $gh$.

Denote by $m$ the midpoint of $gh$.
Take the $(d-2)$-dimensional hyperplane $M$ of $H$ through $m$ orthogonal to $gh$. 

Provide the hyperplane $G$ orthogonal to $gm$ through $g$.
Let $G^+$ be the closed halfspace of $\mathbb{H}^d$ bounded by $G$ which contains $j$, and let $G^-$ be the opposite closed halfspace.

Since $j$ is the unique point of $C$ in the maximum distance from $H$ and by the compactness of $C \cap G^-$ there is a positive number $\kappa^- < \Delta(C)$ for which every point of $C \cap G^-$ is in a distance at most $\kappa^-$ from $H$.
Hence we can rotate $H$ around $M$, in the positive orientation, up to a position $H'$ for which there is positive number $\lambda^-$ strictly between $\kappa^-$ and $\Delta(C)$ such that every point of $C \cap G^-$ is in a distance at most $\lambda^-$ from $H'$.

Clearly, our $H'$ can be simultaneously chosen so close to $H$ that it has empty intersection with $C \cap G^+$, and thus with $C$.
Hence there exists a $\lambda^+$ such that every point of $C \cap G^+$ is in a distance at most $\lambda^+$.

From the two preceding paragraphs we conclude that every point of $C$ is in a distance at most $\lambda= \max \{\lambda^-,\lambda^+\}$ from $H'$.
Since $\lambda < \width_H (C)$, by Proposition 1 applied to $H'$ we conclude that $\width_{H'} (C) < \width_H (C)$.
This contradicts the assumption that $\width_H (C) = \Delta (C)$.
Hence $h$ must be in $C \cap H$.
\end{proof}

\noindent
{\bf Example 2.} In order to find the hight $\eta_x$ of the regular triangle $T$ of side $2x$ we take a vertex $a$ of $T$ and the midpoint $m$ of the opposite side $bc$. 
Look at the triangle $amb$. 
Clearly, at $m$ there is the right angle. 
Thus by the Pythagorean theorem for $\mathbb{H}^2$ we obtain $\ch |ab| = \ch |mb| \cdot \ch |ma|$, this is $\ch 2x = \ch x \cdot \ch \eta_x$. 
Hence $\eta_x = \arch\frac{\ch 2x}{\ch x}$.

\vskip0.3cm
\noindent
{\bf Example 3.} 
Denote by $S_{2x}$ the regular tetrahedron $uvwz$ whose edges are of length $2x$ and by $p$ the projection of $z$ onto the facet $uvw$.
By Lemma 1 of Kellerhals \cite{[Ke]} for $\mathbb{H}^2$ we know that the circumradius of the regular triangle of side $2x$ is $\arsh (\sqrt {4/3} \cdot \sh x)$.
So this is $|up|$. 
Since the triangle $zpu$ has the right angle at $p$, by the Pythagorean theorem for $\mathbb{H}^2$ we obtain 
$\ch 2x = \ch [\arsh (\sqrt \frac{4}{3} \cdot \sh x)] \cdot \ch |zp|$.
This and $|zp| = \width_G (S_{2x})$ for the hyperplane $G$ containing $uvw$ lead to

$$\width_G (S_{2x}) = \arch\frac{\ch 2x}{\ch [\arsh(\sqrt {4/3} \, \sh x)]}. \eqno (1)$$

\section{Bodies of constant width and reduced bodies}

If for every hyperplane supporting a convex body $W \subset \mathbb{H}^d$ the width of $W$ determined by this hyperplane is $\delta$, we say that $W$ is a {\bf body of constant width} $\delta$. 
This definition is analogous to the definitions in $\mathbb{E}^d$ and $\mathbb{S}^d$ (for instance see \cite{[CG]} and \cite{[LaMu]}).
In all the three cases the concept of width of $C$ with respect to a supporting hyperplane (hemisphere in $\mathbb{S}^d$) is analogous. 
By the way, there are different definitions of a body of constant width in $\mathbb{H}^d$. 
Some of them apply different notions of width, as for instance Santal\'o in \cite{[S]}.
Some other take into account a property of a convex body in $\mathbb{E}^d$ equivalent to the definition of a body of constant body there, for instance see the recent paper by B\"or\"oczki and Sagemeister \cite{[BS]}.
A task would be to check the relationships of all such definitions of a body of constant width in $\mathbb{H}^d$, analogously as this is done for $\mathbb{S}^d$ in the last section of \cite{[L22]}.

It is easy to check that balls in $\mathbb{H}^d$  and every Reauleaux polygon in $\mathbb{H}^2$ (defined as in $\mathbb{E}^2$) are bodies of constant width.

From Theorem 1 we obtain the following claim.

\begin{cla}
For any body $W$ of constant width $\delta$ we have $\diam (W) = \delta$.
\end{cla}

\begin{pro} 
Every body of constant width is strictly convex. 
\end{pro}

\begin{proof}
Imagine the opposite to the thesis, i.e., that there exists a body $W$ of constant width $\delta$ which is not strictly convex.
Then a hyperplane $H$ supporting $W$ contains more than one point.
There is a farthest from $H$ ultraparallel hyperplane $J$ supporting $C$.
By Theorem 2 and the sentence just before it there is exactly one point $j \in W \cap J$ and its projection $h$ onto $H$ belongs to $C \cap H$.
Clearly, $|jh| = \delta$.
Recall that $H$ contains a point $t$ different from $h$.
By Lemma 1 we have $|jt| > |jh|$ which implies that $\diam (W) > \delta$.
On the other hand, by Claim 2 we get $\width_H(W) = \diam (W)$.
This contradiction shows that our proposition holds true.
\end{proof}

A few different notions of a body of constant width in $\mathbb{H}^d$ exist.
Some information on them are given in pages 240--242 of the book \cite {[MMO]} by Martini, Montejano and Oliveros.

In analogy to the definition of reduced bodies in Euclidean space $\mathbb{E}^d$ (see \cite{[He]} and the survey article \cite{[LasMar]}) and on the sphere $S^d$ (see \cite{[L15]} and the survey article \cite{[L22]}) we define reduced convex bodies in $\mathbb{H}^d$. 
We say that a convex body $R \subset \mathbb{H}^d$ is {\bf reduced} if $\Delta (Z) < \Delta (R)$ for every convex body $Z \subset R$ different from $R$. 

From Zorn's lemma we see that {\it every convex body $C \subset \mathbb{H}^d$ contains a reduced body whose thickness is $\Delta (C)$}. 

It is easy to show that all regular odd-gons in $\mathbb{H}^2$ are reduced bodies.

\begin{pro}
Every body of constant width in $\mathbb{H}^d$ is reduced.
\end{pro}

\begin{proof} 
Let $W \subset \mathbb{H}^d$ be a body of constant width.
In order to confirm that $W$ is reduced take an arbitrary convex body $Z \subset W$ different from $W$; we intend to show that $\Delta (Z) < \Delta (W)$.

There exists a hyperplane $K$ supporting $Z$ which has non-empty intersection with the interior of $W$.
Take any point $q \in Z$ farthest from $K$, which by Proposition 1 
means that ${\rm dist} (q, K) = {\rm width}_K (Z)$.
Provide an ultraparallel hyperplane $L$ to $K$ which supports $W$ as far as it is possible from $K$ on the opposite side of $K$ than $Z$ is situated. 
By Lemma 1 we have ${\rm dist} (q, K) < {\rm dist} (q, L)$.
By Proposition 1 applied to ${\rm width}_L(W)$ we see that ${\rm width}_L(W) \geq {\rm dist} (q, L)$ (still each most far point of $W$ from $L$ is in a distance at least ${\rm dist} (q, L)$).
Since $W$ is of constant width, every its width is equal to $\Delta (W)$.
In particular, ${\rm width}_L(W) = \Delta (W)$.

By the above observations we obtain $\Delta (W) = {\rm width}_L(W) \geq {\rm dist} (q, L) > {\rm dist} (q, K) = {\rm width}_K(Z) \geq \Delta (Z)$.
From the obtained inequality $\Delta (Z) < \Delta (W)$ we conclude that $W$ is a reduced body.
\end{proof}

Take a reduced body in $\mathbb{H}^2$ with an axis of symmetry and rotate it in  $\mathbb{H}^3$ around this axis. 
It appears that the obtained $3$-dimensional body is reduced (see Fig. 4 for rotating the Reuleaux triangle around its axis of symmetry).
Similarly as in $\mathbb{E}^3$ and $\mathbb{S}^3$.

\vskip0.15cm
\begin{figure}[htbp]        
\hskip0.85cm
\includegraphics[width=5 cm]{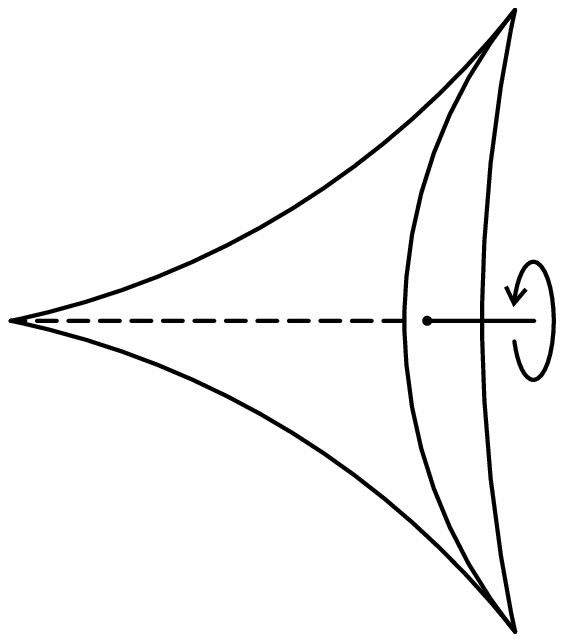} \hskip 1.15cm  
\includegraphics[width=5.62cm]{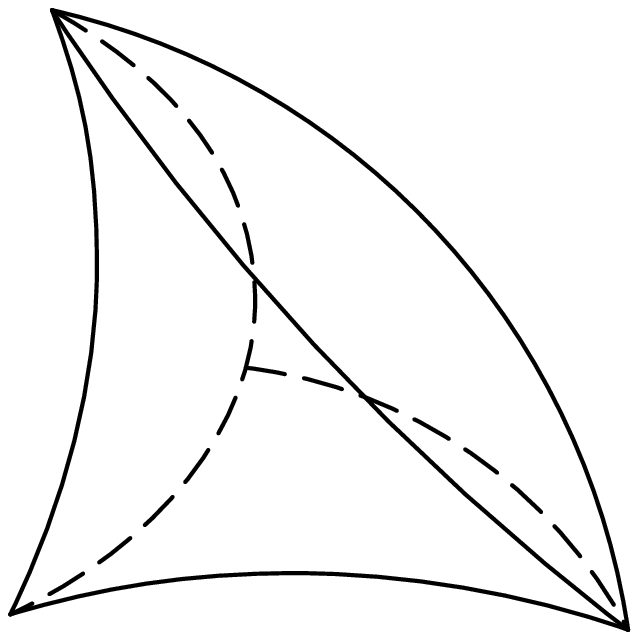} 

\vskip0.2cm
\hskip 0.75cm 
{\bf Fig. 4.} A rotational reduced body \hskip 1.6cm 
{\bf Fig. 5.} A $\frac{1}{8}$-th part of a ball  

\end{figure}

Dissect a ball in $\mathbb{H}^d$ by $d$ orthogonal hyperplanes through its center into $2^d$ closed parts called ${\frac{1}{2^d}}${\bf -th parts of a ball}.
Clearly, the thickness of each of them is equal to the radius of the above ball.
Observe that every ${\frac{1}{2^d}}$-th part of a ball is a reduced body (similarly as in $\mathbb{E}^d$ and $\mathbb{S}^d$, see \cite{[L90]} and \cite{[L15]}).
The $3$-dimensional case is seen in Fig. 5.

\vskip0.2cm
\noindent
{\bf Example 4.} Take $S_{2x}$ from Example 3.
Provide the hyperplane $H$ containing the edge $uv$ and orthogonal to $mn$, where $m$ is the midpoint of $uv$ and $n$ is the midpoint of $wz$ (see Fig. 6). 
By the Pythagorean theorem for the rectangular triangle $nmv$ we obtain $\ch |nv| = \ch|mn| \cdot \ch |mv|$ which means that $\ch|mn| = \frac{\ch |nv|}{\ch |mv|}$, i.e., that $|mn| = \arch \frac{\ch |nv|}{\ch |mv|}$.
Hence $|mn| = \arch \frac{\ch|nv|}{\ch|mv|}$.  
Still $|mv| = x$ and  $\ch |nv| = \ch \eta_x =  \frac{\ch 2x}{\ch x}$ as shown Example 2.
By Proposition 1 the width of $S_{2x}$ determined by $H$ equals to the distance of a farthest point of $S_{2x}$ from $H$.
Such a point is $z$. 
By $h$ denote the projection of $z$ onto $H$.
Then $\width_H(C) = |hz|$. 
Observe that points $m, n, h$ and $z$ are in a hyperplane.
Consider the Lambert quadrilateral $nmhz$.
The right angles are at $n, m$ and $h$.
Thus $\sh |zh| = \sh |mn| \cdot \ch |nz|$.
Hence $|zh| = \arsh( \sh |mn| \cdot \ch |nz|)$.
By Proposition 1 we have $|zh| = \width_H(T)$.
Moreover, taking into account $|nz| = x$ and the established earlier formula for $|mn|$ we obtain 

$$\width_H (S_{2x}) = \arsh[\sh(\arch \frac{\ch 2x}{\ch^2 x})\cdot \ch x]. \eqno (2)$$

\begin{figure}[htbp]        
\hskip2.65cm
\includegraphics[width=9cm]{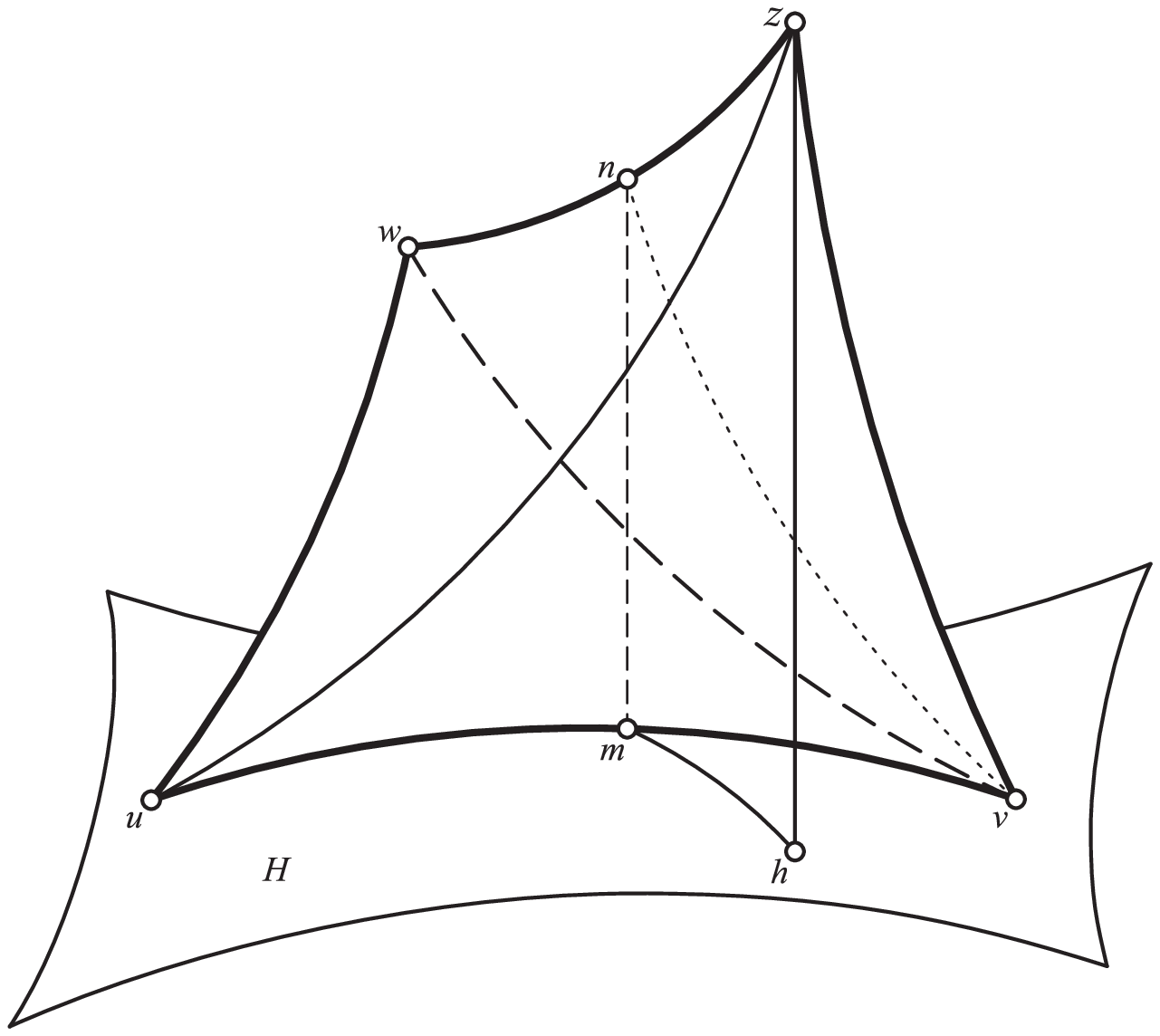} 

\hskip 4.5cm {\bf Fig. 6.} Illustration to Example 4  

\end{figure}

\noindent
{\bf Remark 2.}
In $\mathbb{E}^d$ for $d\geq 3$ every regular simplex is not reduced (see \cite{[MaSw]}). 
The author tried to check the situation in $\mathbb{H}^3$ expecting that sufficiently large tetrahedra are reduced. 
It appears that the opposite is true.
We show this by applying the formulas $\ch (\arch a) = \sqrt{1-a^2}$, $\sh(\arch a) = \sqrt{a^2-1}$, $\ch 2x = 2 \,\ch^2 x -1$ and $\sh x = \sqrt {\ch^2 x -1}$ for the difference of (1) and (2), which after substituting $\lambda = \ch 2x$ leads to 
the inequality $29\lambda^2 - 36\lambda +7 < 0$.
It is never fulfilled for $\lambda>1$, so never for any positive $x$.
So every $S_{2x}$ is not reduced (the author expects that the situation is analogous for every $d > 3$). 
By the way, $\width_H (S_{2x}) / \width_G (S_{2x})$ tends to $1$ as $x \to \infty$.

\vskip0.15cm
We conjecture that for any reduced convex body $R \subset \mathbb{H}^d$ and any supporting hyperplane $H$ such that ${\rm width}_H(C) = \Delta (R)$ there are unique $a \in R \cap H$ and $b \in \bd (R)$ such that $ab \perp H$ and $|ab| = \Delta(R)$. 

Is it true that every smooth reduced body in $\mathbb{H}^d$ is of constant width?
Similarly as in $\mathbb{E}^d$, which is proved by Groemer \cite{[Gr]} (compare Corollary 1 from \cite{[L90]}) and also as in the version for $\mathbb{S}^d$ shown in Theorem 5 from \cite{[L15]}.

\section {Complete and constant diameter bodies}

We say that a convex body $D \subset \mathbb{H}^d$ is of {\bf constant diameter} $\delta$ provided $\diam(D) = \delta$ and for every $p \in \bd(D)$ there is a point $p' \in \bd(D)$ with $|pp'| = \delta$. 
This definition is analogous to the one for $\mathbb{S}^d$ from Section 4 of \cite{[LaMu]}.

Similarly to the traditional notion of a complete set in $\mathbb{E}^d$ (for instance, see \cite{[BF]}, \cite{[CG]} and \cite{[Eg]}) we say that a set $C \subset \mathbb{H}^d$ of diameter $\delta$ is {\bf complete} provided $\diam (C \cup \{x\}) > \delta$ for every $x \not \in C$. 

We omit the proof of the following two claims since they are similar to the proof by Lebesgue \cite{[Le]} in $\mathbb{E}^d$, which is recalled in Section 64 of \cite{[BF]}.

\begin{cla}
Every set of a diameter $\delta$ in $\mathbb{H}^d$ is a subset of a complete set of diameter~$\delta$.
\end{cla}

\begin{cla}
Every complete body $C$ of a diameter $\delta$ in $\mathbb{H}^d$ coincides with the intersection of all balls of radius $\delta$ centered at points of $C$.
\end{cla}

This claim permits to use the term a {\bf complete convex body} for a complete set.

\begin{pro}\label{p'} 
If $C \subset H^d$ is a complete body of diameter $\delta$, then for every $p \in \bd (C)$ there exists $p' \in C$ such that $|pp'|=\delta$. 
\end{pro}

\begin{proof} 
Imagine that the thesis is not true. 
Then $|pp'| < \delta$ for a point $p \in \bd (C)$ and for every point $p' \in C$. 
Since $C$ is compact, there is an $\omega < \delta$ such that  $|pp'| \leq \omega$ for every $p' \in C$. 
Take a ball $B$ of radius $\delta - \omega$ centered at $p$.
Of course, there is a supporting hyperplane $H$ supporting $C$ at $p$.
Observe that the open halfspace bounded by $H$ and containing the interior of $C$ contains the relative interior of the segment $pp'$.
Consequently, the straight line through $p'$ and $p$ contains a point $x \in B$, but out of $C$, in the distance $\delta - \omega$ from $p$.
Thus $|p'x| = \delta$.  
We get $\diam (C \cup \{x\}) \geq \delta$.
Since $x \not \in C$, this contradicts the assumption that $C$ is a complete body of diameter $\delta$.
Hence our proposition holds true.
\end{proof}

This proposition is analogous the Lemma 2 of \cite{[L20]} and holds true also in $E^d$.
Just this proof can be repeated also for every complete body in $E^d$.

\begin{thm}
Let $C \subset \mathbb{H}^d$.
We have $(a) \Longrightarrow (b) \Longleftrightarrow (c)$, where

{\rm (a)} \ $C$ is a body of constant width $\delta$, 
 
{\rm (b)} \ $C$ is a body of constant diameter $\delta$,

{\rm (c)} \ $C$ is complete body of diameter $\delta$.
\end{thm}

\begin{proof}
Show that (a) implies (b).

Assume that $C$ is of constant width $\delta$.
By Theorem 1 we have $\diam (C) = \delta$, which is the first requirement of the definition of a body of constant diameter.

Take an arbitrary point $p \in \bd(C)$.
Let $P$ be any hyperplane supporting $C$ at $p$.
Since $C$ is of constant width $\delta$, we have $\width_P(C) = \delta$. 
Let $p'$ be a farthest point of $C$ from $P$.
By Proposition 1 its distance from $P$ is $\width_P(C)$.
This and $p \in P$ imply $|pp'| \geq \delta$.
Hence by $\diam (C) = \delta$ we get $|pp'| = \width_P (C)$.
We conclude that $p'$ fulfills the second requirement of the definition of a body of constant diameter.

From the above two paragraphs we see that $C$ is a body of constant diameter $\delta$.

\vskip0.1cm
Show that (b) and (c) are equivalent.

Consider a body $C \subset \mathbb{H}^d$ of constant diameter $\delta$.
Take any $r \not \in C$.
Let $B$ be the largest ball whose interior is disjoint with $C$. 
Denote by $\rho$ the radius of this ball and by $p$ the common point of $B$ and $C$.
Since $C$ is of constant diameter $\delta$, there is $p' \in C$ in the distance $\delta$ from $p$.
Clearly, a unique hyperplane supports $B$ and $C$ at $p$.
Hence $\angle rpp' \geq \frac{\pi}{2}$.
Thus by the hyperbolic law of cosines for the triangle $rpp'$ we get 
$|rp'| > |pp'|$ which means that $|rp'| > \delta$. 
Consequently, $\diam(C \cup \{r\}) > \diam (C)$ for every $r \not \in C$. 
So $C$ is complete.

On the other side, if $C \subset \mathbb{H}^d$ is a complete body of diameter $\delta$, then by Proposition 5 we conclude that $C$ is of constant diameter.
\end{proof}

We conjecture that (a), (b) and (c) are equivalent.
In order to confirm this, it is sufficient to show that (b) or (c) imply (a).

\vskip0.2cm
\noindent
Marek Lassak

\noindent
University of Science and Technology

\noindent
85-789 Bydgoszcz, Poland

\noindent
e-mail: lassak@pbs.edu.pl

\end{document}